\documentclass[twocolumn]{autart}    

\usepackage{graphicx}          
                               
\usepackage[numbers]{natbib}        

\usepackage{epsfig} 
\usepackage{amsfonts}
\usepackage{amsmath,amssymb}
\usepackage{comment}
\usepackage{color}
\usepackage{cases}
\usepackage{ascmac}
\usepackage[subrefformat=parens]{subcaption}

\newtheorem{problem}{Problem}

\newtheorem{lemma}{Lemma}
\newtheorem{theorem}{Theorem}
\newtheorem{corollary}{Corollary}
\newtheorem{rmk}{Remark}
\newtheorem{example}{Example}

\newcommand{\R}{{\mathbb R}}

\newcommand{\calF}{{\mathcal F}}

\newcommand{\calL}{{\mathcal L}}

\newcommand{\calS}{{\mathcal S}}    
    
\newcommand{\calU}{{\mathcal U}}

\newcommand{\bbE}{{\mathbb E}}

\newcommand{\bbP}{{\mathbb P}}    
    
\newcommand{\bbR}{{\mathbb R}}

\newcommand{\bbU}{{\mathbb U}}

\newcommand{\argmax}{\mathop{\rm arg~max~}\limits}

\newcommand{\ito}[1]{{\color{black}{#1}}} 
\newcommand{\ike}[1]{{\color{black}{#1}}} 
\newcommand{\ikeda}[1]{{\color{black}{#1}}} 
\newcommand{\rev}[1]{{\color{black}{#1}}}

\newcommand{\modi}[1]{{\color{black}{#1}}} 
\newcommand{\state}[1]{{{#1}}} 
\newcommand{\mnr}[1]{{\color{black}{#1}}}

\newcommand{\tx}{{t,x}}

\begin{document}

\begin{frontmatter}

\title{Sparse Optimal Stochastic Control\thanksref{footnoteinfo}
\vspace{-1.0cm}} 

\thanks[footnoteinfo]{This paper was not presented at any IFAC 
meeting. Corresponding author K.~Kashima. Tel. +81-75-753-5512. \\
©~2021. This manuscript version is made available under the CC-BY-NC-ND 4.0 license https://creativecommons.org/licenses/by-nc-nd/4.0/
}

\author[Kyoto]{Kaito Ito}\ead{ito.kaito@bode.amp.i.kyoto-u.ac.jp},    
\author[Kyushu]{Takuya Ikeda}\ead{t-ikeda@kitakyu-u.ac.jp},               
\author[Kyoto]{Kenji Kashima}\ead{kk@i.kyoto-u.ac.jp}  

\address[Kyoto]{Graduate School of Informatics, Kyoto University, 
	Kyoto, Japan}  
\address[Kyushu]{Faculty of Environmental Engineering, The University of Kitakyushu, 
	Kitakyushu, Japan}             

\begin{keyword}                           
sparsity, non-smooth optimal control, bang-off-bang control, dynamic programming, viscosity solution           
\end{keyword}                             

\begin{abstract}                          
In this paper, we investigate a sparse optimal control of continuous-time \ito{stochastic} systems.
We adopt the dynamic programming approach
and analyze the optimal control via the value function.
Due to the non-smoothness of the $L^0$ cost functional, in general, 
the value function is not differentiable in the domain.
Then, we characterize the value function as a viscosity solution to the associated Hamilton-Jacobi-Bellman (HJB) equation.
Based on the result, 
we derive a necessary and sufficient condition for the $L^0$ optimality, 
which immediately gives the optimal feedback map.
\modi{Especially for control-affine systems}, we consider the relationship with $L^1$ optimal control problem
and show an equivalence theorem.
\vspace{-.5cm}
\end{abstract}

\end{frontmatter}

\section{Introduction}
This work investigates an optimal control problem for non-linear \ike{stochastic} systems 
with the $L^0$ control cost.
This cost functional penalizes the length of the support of control variables,
and the optimization based on the criteria tends to make the control input identically zero on a set with positive measures.
\rev{Consequently,} the optimal control is switched off completely on parts of the time domain.
Hence, this type of control is also referred to as 
{\em sparse optimal control}.
For example, this optimal control framework is applied to 
actuator placements~\cite[]{Sta09,HerStaWac12,KunPieVex14},
networked control systems~\rev{\cite[]{NagQueOst14,IkeKas18_node,Olshevsky2020}},
and \mnr{discrete-valued control}~\cite[]{Ikeda2016}, 
to name a few.
The sparse optimal control involves the discontinuous and non-convex cost functional.
Then, in order to deal with the difficulty of analysis,
some relaxed problems with the $L^p$ cost functional have been often investigated,
akin to methods used in compressed sensing applications~\cite[]{Don06}.

{\textit{{Literature review:}}}\ \
\ito{For {\em deterministic} control-affine systems}, the $L^1$ cost functional is analyzed with an aim to show the relationship between the $L^0$ optimality and the $L^1$ optimality, and an equivalence theorem is derived in~\cite[]{NagQueNes16}.
In~\cite[]{IkeKas18}, the result is extended to \ito{deterministic} general linear systems including infinite-dimensional systems.
The $L^1$ control cost is also considered in~\cite[]{AltSch15,Vos06,AthFal}.
In~\cite[]{IkeKas18}, the sparsity properties of optimal controls for the $L^p$ cost with $p \in (0,1)$ is discussed.
The authors investigated this problem from a dynamic programming viewpoint~\cite[]{IkeKas_ECC19}. 
When it comes to {\em stochastic} systems, in~\cite[]{ExaThe18}, a finite horizon optimal control problem with the $L^1$ cost functional for stochastic systems is dealt with and \ito{the authors propose} sampling-based algorithm to solve the problem utilizing forward and backward stochastic differential equations. However, it is not obvious that the $L^1$ optimal control achieves the desired sparsity.
To the best of the authors' knowledge, our preliminary work~\rev{\cite[]{Ito2020ifac}} (Theorem~\ref{thm:conti_V_s} below) on the continuity of the value function is the only theoretical result on $L^0$ optimal control of stochastic systems.

{\textit{{Contribution:}}}\ \
The goal of this work is to obtain the \rev{sparse optimal} feedback \rev{map} (Theorem~\ref{thm:discrete}), where the optimal control input has the bang-off-bang property and to reveal the equivalence between the $L^0$ optimality and the $L^1$ optimality for \modi{control-affine} stochastic systems (Theorem~\ref{cor:L0_L1_s}).
To this end, we utilize the dynamic programming. 
In the present paper, we first characterize our value function as a \emph{viscosity solution} to the Hamilton-Jacobi-Bellman (HJB) equation~\cite[]{Fleming2006,Yong1999}.
Based on the result,
we show a \modi{necessary and sufficient} condition for the $L^0$ optimality \modi{(Theorem~\ref{thm:feedback})},
which immediately gives an optimal feedback \rev{map}.
In addition, 
a sufficient condition for the value function to be a classical solution to the HJB equation \rev{(i.e., a solution that satisfies the HJB equation in the usual sense)} is given via the equivalence, \modi{while, in general, for the deterministic case, we cannot ensure the differentiability of the value function.}

\modi{In the stochastic case, the HJB equation becomes a second-order equation compared to that of the deterministic case, 
	and hence the results for the deterministic systems~\cite{IkeKas_ECC19} cannot be directly applied to the stochastic case.
	Indeed, several difficulties arise due to the stochasticity. For example, the analysis of the deterministic $L^0$ optimality of a control in~\cite{IkeKas_ECC19} heavily relies on the local Lipschitz continuity of the value function, which means the almost everywhere differentiability. On the other hand, the value function for the stochastic $L^0$ optimal control is at most locally 1/2-H\"{o}lder continuous~\cite[]{Ito2020ifac}. This implies that we need a quite different approach. Even for the problem formulation, we must be careful about the probability space we work on to correctly apply the dynamic programming principle.
}

In order to demonstrate the practical usefulness of \mnr{our theoretical results}, an example is exhibited; see Example \ref{ex:2} for more details.

\begin{figure}[tb]
			\centering
			\includegraphics[scale=0.4]{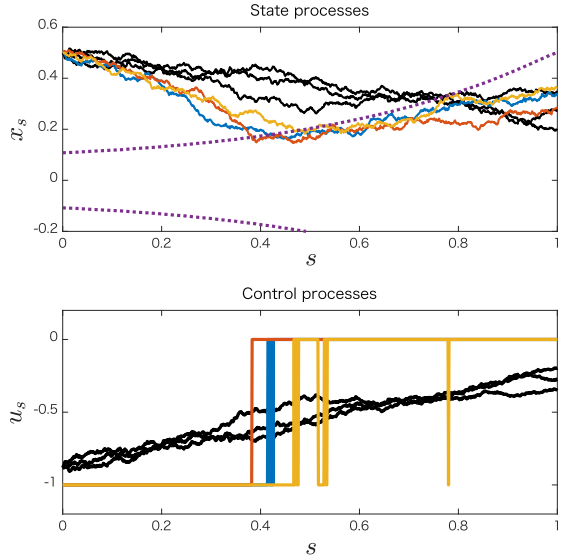}
			\caption{The colored lines except black are the sample paths of the $L^0$ optimal state process (top, solid) and control process (bottom, solid), and the switching \rev{boundary} (top, dotted). The same color indicates the correspondence between the sample paths of the state process and the control process.
				The black lines are the sample paths of the $L^2$ optimal state process (top) and control process (bottom).}
			\label{fig:sample_path}
		\end{figure}

\begin{example}
Consider the following stochastic system:
\begin{equation}\label{eq:simple_s}
d\state{x}_s = c\state{x}_s ds + u_s ds + \sigma dw_s, \quad 0\leq s \leq T
\end{equation}
where $\{\state{x}_s\}$ is a real-valued state process,  $\{u_s\}$ is a control process, and $\{w_s\}$ is a Wiener process. We take $c=1,\ \sigma=0.1,\ T=1$, and $\state{x}_0 = 0.5$. 
Then, the black lines in Fig.~\ref{fig:sample_path} show sample paths of the optimal control input and corresponding state trajectories that minimize $\bbE\left[ \int_0^1 |u_s|^2 ds + \state{x}_1^2 \right]$. It is well known that this minimum energy  control is given by linear state-feedback control, and hence it takes non-zero values almost everywhere.
On the contrary, our problem can deal with
the sparse optimal control that minimizes $\bbE\left[ \int_0^1 |u_s|^0 ds + \state{x}_1^2 \right]$
with the constraint $|u_s|\le 1$ where $0^0 = 0$. The first term represents the length of time that the control takes non-zero values. \emph{\rev{Theorem~\ref{thm:discrete}} reveals that the optimal control input takes only three values of $\{-1,0,1\}$, and enables us to numerically compute the state-feedback \rev{map} from $\state{x}_s$ to $u_s\in \{-1,0,1\}$}.
The colored lines show the result of $L^0$ optimal control, whose input trajectories are sparse while the variance of the state is small enough. Note that the purple dotted lines show the boundary of the bang-off-bang regions.
\hfill$\lhd$
\end{example}

{\textit{{Organization:}}}\ \
The remainder of this paper is organized as follows.
In Section~\ref{sec:math}, we give mathematical preliminaries for our subsequent discussion.
\ike{In Section~\ref{sec:formulation}, we describe the system model and formulate the sparse optimal control problem for stochastic systems.}
\ike{Section~\ref{sec:preliminary_analysis} is devoted to the general analysis of the stochastic optimal control with the discontinuous $L^0$ cost.}
We first characterize the value function as a viscosity solution to the associated HJB equation and next show a \mnr{necessary and sufficient} condition for the $L^0$ optimality.
\ike{Section~\ref{sec:char_stochastic} characterizes the \rev{sparse optimal stochastic} control.
	We show the relationship with the $L^1$ optimization problem 
	and some basic properties of the \rev{sparse optimal stochastic} control 
	for control-affine systems with box constraints.
}
In Section~\ref{sec:conclusion} we offer concluding remarks.

\section{Mathematical preliminaries}
\label{sec:math}
This section reviews notation that will be used throughout the paper.

Let $N$, $N_1$, and $N_2$ be positive integers.
For a matrix $M\in{\mathbb{R}}^{N_1\times N_2}$, 
$M^{\top}$ denotes the transpose of $M$.
\ike{For a matrix $M\in{\mathbb{R}}^{N\times N}$, 
	$\mathrm{tr}(M)$ denotes the trace of $M$. 
	Denote by $\calS^N$ the set of all symmetric $N\times N$ matrices
	and by $\calS_+^N$ the set of all positive semidefinite matrices.
	Denote the Frobenius norm of $M\in\R^{N_1\times N_2}$ by $\|M\|$, i.e., $\|M\| \triangleq \sqrt{\textrm{tr}(M^\top M)}$.}
For a vector \mnr{$a=[a^{(1)}, a^{(2)}, \dots, a^{(N)}]^{\top}\in\mathbb{R}^{N}$},
we denote the Euclidean norm by 
$\|a\|\triangleq (\sum_{i=1}^{N} (a^{(i)})^{2})^{1/2}$ and
the open ball with center at $a$ and radius $r>0$ by $B(a,r)$, i.e.,
$B(a,r)\triangleq\{x\in\mathbb{R}^N: \|x-a\| < r\}$.
We denote the inner product of $a\in\mathbb{R}^N$ and $b\in\mathbb{R}^N$ by $a\cdot b$.

\mnr{For $p\in\{0,1\}$ and a continuous-time signal $u_s=[u_s^{(1)}, u_s^{(2)}, \dots, u_s^{(N)}]^{\top}\in{\mathbb{R}}^N$ over a time interval $[t, T]$, the {\em $L^p$ norm} of $u = \{u_s\}_{t\le s \le T}$ is defined by
\begin{align*}
&\|u\|_{0} \triangleq \sum_{j=1}^{N}\mu_{L}(\{s\in[t, T]: u_s^{(j)}\neq0\}),\\
&\|u\|_{1} \triangleq \sum_{j=1}^{N} \int_{t}^{T} |u_s^{(j)}| ds,
\end{align*}
}with the Lebesgue measure $\mu_{L}$ on $\mathbb{R}$.
The $L^0$ norm is also expressed by
$\|u\|_{0}=\int_{t}^{T}\psi_0(u_s) ds$,
where $\psi_0:\mathbb{R}^N\to\mathbb{R}$ is a function 
that returns the number of non-zero components, i.e.,
\[
\psi_0(a) \triangleq \sum_{j=1}^{N} \mnr{|a^{(j)}|^0}, \ a\in \bbR^N
\]
with $0^0=0$.

For a given set $\Omega\subset\mathbb{R}^N$,
$C(\Omega)$ denotes the set of all continuous functions on $\Omega$.
\mnr{For $T > 0$, $C^{1,2} ((0,T)\times \bbR^N)$ denotes the set of all functions $\phi$ on $(0,T)\times \bbR^N$ whose partial derivatives $\frac{\partial \phi}{\partial s}, \frac{\partial \phi}{\partial x^{(i)}}, \frac{\partial^2 \phi}{\partial x^{(i)} \partial x^{(j)}},  i,j =1,\ldots, N,$ exist and are continuous on $(0,T)\times \bbR^N$.
Denote by $C^{1,2} ([0,T) \times \bbR^N)$ the set of all $\phi \in  C^{1,2} ((0,T)\times \bbR^N) \cap C([0,T)\times \bbR^N)$ such that $\frac{\partial \phi}{\partial s}, \frac{\partial \phi}{\partial x^{(i)}}, \frac{\partial^2 \phi}{\partial x^{(i)} \partial x^{(j)}},  i,j =1,\ldots, N,$ can be extended to continuous functions on $C([0,T) \times \bbR^N)$.
For $\phi\in C^{1,2} ([0,T)\times \bbR^N)$, $\phi_t$ denotes the partial derivative with respect to the first variable, $D_x \phi$ denotes the gradient with respect to the last $N$ variables, and $D_x^2 \phi$ denotes the Hessian matrix with respect to the last $N$ variables.
}
\mnr{For $p \ge 2$, denote by $C_p^{1,2} ([0,T]\times \bbR^N)$ the set of all $\phi \in C^{1,2}([0,T) \times \bbR^N) \cap C([0,T]\times \bbR^N)$ satisfying
	\begin{equation}\label{ineq:phi_dphi_growth}
	\|\rho(\tx)\| \le K(1+\|x\|^p), \  (\tx)\in [0,T]\times \bbR^N 
	\end{equation}
	for some constant $K > 0$ and any $\rho \in \{\phi, D_x \phi, D_x^2 \phi, \phi_t \}$.
A function $\rho : [0,T]\times \bbR^N \rightarrow \bbR$ is said to satisfy a polynomial growth condition or to be at most polynomially growing if there exist constants $K >0$ and $p \ge 2$ such that \eqref{ineq:phi_dphi_growth} holds.
}

\ito{Let $\alpha \in (0,1]$. A function $f:\bbR^{N_1} \rightarrow \bbR^{N_2}$ is called $\alpha$-H\"{o}lder continuous if there exists a constant $L>0$ such that, for all $x,y\in \bbR^{N_1}$,
	$\|f(x) - f(y)\| \le L\|x-y\|^\alpha$.
	Especially when $\alpha = 1$, $f$ is called Lipschitz continuous. A function $f$ is called locally $\alpha$-H\"{o}lder continuous if for any $x\in \bbR^{N_1}$, there exists a neighborhood $U_x$ of $x$ such that $f$ restricted to $U_x$ is $\alpha$-H\"{o}lder continuous.}

\modi{The notation $o(s)$ denotes a real-valued function $f$ defined on some subset of $\bbR$ such that $\lim_{s\rightarrow 0} f(s)/s = 0$. }

\ito{For $0\le t \le T$, let $(\Omega, \calF, \{\calF_s\}_{s\ge t}, \bbP)$ be a filtered probability space, and $\bbE$ be the expectation with respect to $\bbP$. For $S = \bbR^N {\rm ~or~} \calS^N$, denote by \mnr{$\calL_{\calF}^2 (t,T;S)$} the set of all $\{\calF_s\}_{s\ge t}$-adapted $S$-valued processes $\{X_s\}_{s\ge t}$ such that $	\bbE\left[ \int_t^T \|X_s\|^{\mnr{2}} ds \right] < +\infty$. 
	In what follows, we omit the subscript of stochastic processes when no confusion occurs, e.g., $\{X_s\} = \{X_s\}_{s \ge t}$.

}

\section{Problem formulation}
\label{sec:formulation}
\ike{This paper considers the sparse optimal control for stochastic systems.
	This section provides the system description and formulates the main problem.}

\ike{We consider the following stochastic system where the state is governed by a stochastic differential equation valued in $\bbR^n$:
	\begin{align}\label{eq:stochastic_system}
	\begin{split}
	&d\state{x}_s = f(\state{x}_s, u_s) ds +  \sigma(\state{x}_s, u_s) dw_s, \quad s>t, \\
	&\state{x}_t = x.
	\end{split}
	\end{align}
	\modi{The initial value $x\in \bbR^n$ is deterministic, and $\{w_s\}$ is a $d$-dimensional Wiener process.} The range of the control $\bbU \subset \bbR^m$ is a compact set that contains $0\in \bbR^m$, 
	and we fix a finite horizon $0<T<\infty$.

	We are interested in the optimal control that minimizes the cost functional
	\begin{equation}\label{eq:cost_stoc}
	J^{\sf s} (t,x,u) \triangleq \bbE \left[ \int_t^T \psi_0 (u_s) ds + g(\state{x}_T) \right].
	\end{equation}
	
	\modi{We assume the following conditions for functions $f, \sigma, g$:
	\begin{enumerate}
		\item[$(A_1)$]
		The functions $f$ and $\sigma$ are globally Lipschitz, namely, there exist positive constants $L$, $\bar{M}$ and a nondecreasing function $\bar{m} \in C([0,+\infty))$ such that $f:\bbR^n \times \bbU \to \bbR^n$ and $\sigma:\bbR^n \times \bbU \to \bbR^{n \times d}$ satisfy the following condition:
			\begin{align}
			&\|f(\state{x},u) - f(y,v)\|  + \|\sigma(\state{x},u) - \sigma(y,v)\|   \nonumber\\
			&\quad \leq L \|\state{x}-y\| + \bar{m}(\|u-v\|)
			\label{eq:Lip_fsigma}
			\end{align}
			for all $\state{x},y \in \bbR^n$, $u, v \in \bbU$,
			where $\bar{m}(\cdot) \le \bar{M}$ and $\bar{m}(0) = 0$; \label{ass:fsigma_Lip}  
		\item[$(A_2)$]
		There exist constants $\hat{C}>0$ and $p\ge 2$ such that $g:\bbR^n \to \bbR$ satisfies the following growth condition:
		\begin{equation}\label{ineq:g_growth}
		|g(\state{x})| \le \hat{C}(1+\|\state{x}\|^p)
		\end{equation}
		for all $\state{x} \in \bbR^n$;
		\label{ass:g_growth}
		\item[$(A_3)$]
		$g:\bbR^n \to \bbR$ is continuous.
		\label{ass:g_cont_stoc}
	\end{enumerate}
	}
	Given a probability space with the filtration $\{\calF_s\}_{s\ge t}$ generated by a Wiener process, \modi{Assumption~$(A_1)$ ensures the existence and uniqueness of a strong solution to the stochastic differential equation~\eqref{eq:stochastic_system} 
	with any initial condition $\state{x}_t = x, (t,x)\in [0,T] \times \bbR^n$, and any $\{\calF_s\}_{s\ge t}$-progressively measurable and $\bbU$-valued control process $\{u_s\}$.}
	In addition, under assumptions $(A_1)$ and $(A_2)$, 
	the cost functional $J^{\sf s} (t,x,u)$ is finite; see Appendix~\ref{app:moment}.
	Assumption~$(A_3)$ is introduced to show the continuity of the value function defined later in \eqref{eq:value_func}.

	\modi{
	    For our analysis, we utilize the method of dynamic programming.
		\modi{In order to establish the dynamic programming principle (Lemma~\ref{lem:DP_s}), we need to consider a family of optimal control problems with different initial times and states $(t,x)\in [0,T] \times \bbR^n$ along a state trajectory. Let us consider a state trajectory starting from $x_0 = x$ on a filtered probability space $(\Omega, \calF, \{\calF_s\}_{s\ge 0}, \bbP)$. For any $s > 0$, $x_s$ is a random variable. However, an $\{\calF_s\}_{s\ge 0}$-progressively measurable control $\{u_s\}$ knows the information of the system up to the current time. In particular, the current state $x_s$ is deterministic under a {\em different} probability measure $\bbP(\cdot | \calF_s)$. This observation naturally leads us to vary the probability spaces as well as control processes}; for details see e.g.,~\cite{Yong1999,Nisio2014,Fabbri2017}. 
		For this reason, we adopt the so-called {\em weak} formulation of the stochastic optimal control problem; see also Remark~\ref{rmk:weak_formulation}. 
	}

	\ito{For each fixed $t\in [0,T)$, we denote by $\calU^{\sf s}[t,T]$ the set of all 5-tuples $(\Omega, \calF, \bbP, \{w_s\}, \{u_s\})$ satisfying the following conditions:
		\begin{itemize}
			\item[(i)] $(\Omega, \calF, \bbP)$ is a complete probability space,
			\item[(ii)] $\{w_s\}$ is \ito{a} $d$-dimensional Wiener process on $(\Omega, \calF, \bbP)$ over $[t,T]$ (with $w_t = 0$ almost surely),
			\item[(iii)] The control $\{u_s\}$ is an \modi{$\{\calF_s\}_{s\ge t}$}-progressively measurable  and $\bbU$-valued process on $(\Omega, \calF, \bbP)$ where $\calF_s$ is the $\sigma$-field generated by $\{w_r\}_{t\le r \le s}$.
		\end{itemize}
		
		For $(\Omega, \calF, \bbP, \{w_s\}, \{u_s\}) \in \calU^{\sf s}[t,T]$, we call $\{u_s\}$ and \mnr{$(\Omega, \calF, \bbP, \{w_s\})$ an {\em admissible} control process and a \modi{reference probability space}, respectively.}
		For notational simplicity, we \mnr{sometimes} write $u \in \calU^{\sf s}[t,T]$ instead of $(\Omega, \calF, \bbP, \{w_s\}, \{u_s\}) \in \calU^{\sf s}[t,T]$. Note that in $\eqref{eq:cost_stoc}$ the expectation $\bbE$ is with respect to $\bbP$.
		For given $(t,x)\in [0,T]\times \bbR^n$ and $u \in \calU^{\sf s}[t,T]$, we denote by $\{\state{x}_s^{t,x,u}\}_{t\le s \le T}$ the \mnr{unique} solution of \eqref{eq:stochastic_system}. \mnr{When there is no confusion, we omit the subscript.}}
	
	Then, we are ready to formulate the main problem as follows:
	\begin{problem}
		\label{prob:main_stochastic}
		Given $x\in\mathbb{R}^{n}$, $T>0$, and $t\in[0,T]$, 
		find \ito{a 5-tuple $u\in \calU^{\sf s}[t,T]$} that solves
		\begin{equation*}
		\begin{aligned}
		& \underset{u}{\text{minimize}}
		& & J^{\sf s}(t,x,u)\\
		& \text{subject to}
		& & d\state{x}_s = f(\state{x}_s, u_s) ds  + \sigma(\state{x}_s,u_s) dw_s,\\
		& & & \state{x}_t=x,\\
		& & & u\in \mathcal{U^{\sf s}}[t,T].
		\end{aligned}
		\end{equation*}
		\hfill$\lhd$
	\end{problem}

	The value function for Problem~\ref{prob:main_stochastic} is defined by
	\begin{equation}\label{eq:value_func}
	V^{\sf s}(t,x) \triangleq \inf_{u\in \calU^{\sf s}[t,T]} J^{\sf s}(t,x,u),\ (t,x)\in [0,T]\times \bbR^n .
	\end{equation}

}

\ito{
	\begin{rmk}\label{rmk:weak_formulation}
		In Problem \ref{prob:main_stochastic}, we vary probability spaces. This problem formulation is called a {\em weak formulation}. 
		On the other hand, the problem where we fix a probability space for each initial time and state $(t,x)\in [0,T]\times \bbR^n$ and vary only control processes is referred to as a {\em strong formulation}, \modi{which is natural from the practical point of view. Despite the difference in the settings}, it is known that, under some conditions, the value function of the weak formulation \ito{coincides} with the one of the strong formulation; see~\cite[]{Fleming2006}. \modi{In this paper, under some assumptions, we will show that, for any given reference probability space, we can design an optimal state-feedback controller in Corollary~\ref{cor:state-feedback}. This result bridges the gap between the weak formulation and the strong formulation.}
		\modi{Lastly, we would like to emphasize that the term ``weak" refers only to the fact that the probability spaces vary and not to the concept of solution of the stochastic differential equation \eqref{eq:stochastic_system}. In fact, once we fix $u\in \calU^{\sf s}[t,T]$, then the solution is defined on the same probability space.}
		\hfill$\lhd$
\end{rmk}}


\section{General analysis of stochastic optimal control with discontinuous input cost \rev{functional}}\label{sec:preliminary_analysis}
This section is devoted to the preliminary analysis of the stochastic $L^0$ optimal control problem. We first characterize the value function as a viscosity solution to the associated HJB equation. Then, we derive \mnr{a necessary and sufficient} condition for the $L^0$ optimality.

\subsection{Characterization of the value function}
In what follows, we show that the value function $V^{\sf s}$ is a viscosity solution to the associated HJB equation. \modi{The definition of a viscosity solution appears in Appendix~\ref{app:viscosity}.}
The HJB equation~\modi{\cite{Yong1999}} corresponding to the stochastic system \modi{\eqref{eq:stochastic_system}} is given by
\begin{numcases}{}
-v_t(t,x) + H^{\sf s}(x, D_x v(t,x), D_x^2 v(t,x)) = 0, \label{eq:SHJB} \\ \hspace{4cm} (t,x) \in [0,T)\times \mathbb{R}^n, \nonumber \\
v(T,x) = g(x), \quad x\in \mathbb{R}^n, \label{eq:SHJB_initial}
\end{numcases}
where $H^{\sf s}:\mathbb{R}^n\times\mathbb{R}^n \times \calS^n \to\mathbb{R}$ is defined by
\begin{align}\label{eq:H_s}
\begin{split}
H^{\sf s}(x, p,M) &\triangleq\underset{u\in \mathbb{U}}{\sup}\Bigl\{-f(x,u) \cdot p \\
&\quad\quad - \frac{1}{2} \mathrm{tr}(\sigma \sigma^\top (x,u) M)  - \psi_0(u)\Bigr\}.
\end{split}
\end{align}
We first \rev{introduce the result for} the continuity of the value function~\rev{\cite[]{Ito2020ifac}}. \modi{The main difficulty in the analysis is that the state of the system \eqref{eq:stochastic_system} is unbounded due to the stochastic noise.}
\begin{theorem}\label{thm:conti_V_s}
	Fix $T>0$.
	Under assumptions $(A_1),(A_2)$, and $(A_3)$, the value function $V^{\sf s}$ \modi{defined by \eqref{eq:value_func}}, is continuous on $\mnr{[0,T]\times \bbR^n}$. If in addition the terminal cost $g$ is Lipschitz continuous, then $V^{\sf s}(\mnr{t,x})$ is Lipschitz continuous in $x$ uniformly in $t$, and locally $1/2$-H\"{o}lder continuous in $t$ for each $x$.
	\hfill$\lhd$
\end{theorem}

\begin{rmk}\label{rmk:conti_Vs}
	\ito{Note that
		the Lipschitz continuity of $g$ shows the local Lipschitz continuity of the value function for deterministic systems~\cite[Theorem 1]{IkeKas_ECC19}, which ensures that the value function is differentiable almost everywhere. On the other hand, we cannot expect the local Lipschitz continuity of the value function $V^{\sf s}$ in the stochastic case even under the Lipschitz continuity of $g$. This is essentially because $\int_0^t \sigma dw$ is only of order $t^{1/2}$.
		\hfill$\lhd$
	}	
\end{rmk}

The dynamic programming principle plays an important role in proving that the value function is a viscosity solution to the HJB equation. \modi{Since the proof is similar to \cite[Chapter 4, Theorem 3.3]{Yong1999}, it is omitted.}
\begin{lemma}\label{lem:DP_s}
	Fix \modi{any $T>0$ and any $\tau \in [0,T]$}. Assume $(A_1)$ and $(A_2)$. Then, \modi{the value function \eqref{eq:value_func} satisfies}
	\[ V^{\sf s}(\mnr{t,x}) = \inf_{u\in \calU^{\sf s}[t,T]} \bbE\left[ \int_t^{\tau} \psi_0(u_s) ds + V^{\sf s}(\tau, \state{x}_{\tau}^{t,x,u}) \right] \]
	for all $(t,x)  \in [0, \tau]\times \bbR^n$.
	\hfill$\lhd$
\end{lemma}
\ito{\mnr{According to the definition of a viscosity solution in Appendix~\ref{app:viscosity}, we have to check the inequalities \eqref{ineq:subsolution} and \eqref{ineq:supersolution} for {\em any} smooth function $\phi$. However, this requirement is too strong for our analysis.
Fortunately, it is possible to restrict the class of $\phi$ to be considered by the following lemma.}

	\begin{lemma}\label{lem:value_growth}
		Assume $(A_1),(A_2)$, and $(A_3)$. Then, \modi{the value function \eqref{eq:value_func}} satisfies the polynomial growth condition, i.e., for some constant $\hat{C}_p > 0$,
		\begin{equation}\label{ineq:value_growth}
		|V^{\sf s}(\mnr{t,x})| \le \hat{C}_p (1 + \|x\|^p), \ (\tx) \in [0,T]\times \bbR^n
		\end{equation}
		holds where $p\ge 2$ satisfies \eqref{ineq:g_growth}. 
		\modi{In addition, if \eqref{ineq:subsolution} and \eqref{ineq:supersolution} where $v$ and $H$ are replaced by $V^{\sf s}$ and $H^{\sf s}$, respectively, are satisfied for any $\phi\in C_p^{1,2} ([0,T] \times \bbR^n)$, then $V^{\sf s}$ is a viscosity solution to the HJB equation \eqref{eq:SHJB} with a terminal condition \eqref{eq:SHJB_initial}.}
	\end{lemma}
	\begin{pf}
		First, we derive the polynomial growth condition of $V^{\sf s}$.
		By \mnr{Assumption}~$(A_2)$,
		\begin{align} 
		|V^{\sf s} (\tx)| &\le \bbE[|g(\bar{\state{x}}_T)|] \nonumber\\
		&\le \bbE[ \hat{C}(1+\|\bar{\state{x}}_T\|^p)] \label{ineq:upper_bound_valuefunc}
		\end{align}
		holds, where $\hat{C}>0$ and $p\ge 2$ are constants that satisfy \eqref{ineq:g_growth}, and
		$\{\bar{\state{x}}_s\}$ is the solution of the uncontrolled system:
		\[ d\bar{\state{x}}_s = f(\bar{\state{x}}_s, 0)ds  + \sigma(\bar{\state{x}}_s, 0) dw_s, \quad \bar{\state{x}}_t = x . \]
		\modi{Combining the inequality \eqref{ineq:upper_bound_valuefunc} and inequality \eqref{ineq:moment_upper} of Lemma \ref{lem:estimate_moment} in Appendix~\ref{app:moment}}, we obtain \eqref{ineq:value_growth}.
		
		Next, \modi{note that by the definition \eqref{eq:value_func} of the value function $V^{\sf s}$, it satisfies the terminal condition \eqref{eq:SHJB_initial}.}
		Moreover, thanks to the continuity of $V^{\sf s}$ and the derived growth condition \eqref{ineq:value_growth}, we can apply \modi{Theorem 3.1} of \cite{Nisio2014}, that is,
		\modi{$V^{\sf s}$ is a viscosity subsolution (resp. supersolution) of \eqref{eq:SHJB} if \eqref{ineq:subsolution} (resp. \eqref{ineq:supersolution}) holds for any $\phi \in C^{1,2} ([0,T)\times \bbR^n) \cap C([0,T]\times \bbR^n)$ satisfying, for some $R>0$,
		\begin{equation}\label{eq:ass_pgrowth}
		\phi(\tx) = c_p (1+\|x\|^p) ~~ {\rm for} ~ t\in [0,T] , \ \|x\| \ge  R  ,
		\end{equation}
		where $c_p = \hat{C}_p$ (resp. $c_p = -\hat{C}_p$).}
		This implies that for some large $K^{(1)}>0$,
		\[
		|\phi(\tx)| \le K^{(1)} (1+\|x\|^p) , \  (\tx) \in \mnr{[0,T]}\times \bbR^n ,
		\]
		noting that $\phi$ is continuous, \modi{and therefore $|\phi|$ attains a maximum on $[0,T] \times \{x\in \bbR^n: \|x\|\le R\}$}. Moreover, $\eqref{eq:ass_pgrowth}$ gives
		\[
 			D_x \phi (\tx) = p c_p \|x\|^{p-2} x ~~ {\rm for} ~ t\in \mnr{[0,T]} , \ \|x\| \ge  R,  
		\]
		and hence, for some constant $K^{(2)}>0$, it holds that
		\[ \| D_x \phi(\tx) \|  \le K^{(2)} (1+\|x\|^p), ~ (\tx) \in \mnr{[0,T]} \times \bbR^n  , \]
		noting that $D_x \phi$ is continuous.
		Likewise, for some constants $K^{(3)},K^{(4)} >0$, it holds that
		\begin{align*}
		&\|D_x^2 \phi(\tx) \| \le K^{(3)} (1+\|x\|^p), \\
		 &|\phi_t (\tx)| \le K^{(4)} (1+\|x\|^p),~~  (\tx) \in \mnr{[0,T]} \times \bbR^n .
		\end{align*}
		This completes the proof.
		\hfill$\Box$
	\end{pf}

}

Then, we are ready to prove that our value function is a viscosity solution to the associated HJB equation.
\begin{theorem}\label{thm:stochastic_viscosity}
	Fix $T>0$. Assume $(A_1), (A_2)$, and $(A_3)$.
	Then, the value function \modi{\eqref{eq:value_func}} is a viscosity solution to the HJB equation \eqref{eq:SHJB} with a terminal condition \eqref{eq:SHJB_initial}.
\end{theorem}
\begin{pf}
	Note that \mnr{$H^{\sf s}$ is continuous (see Lemma~\ref{lem:continuity_Hs} in Appendix~\ref{app:continuity_H}), 
	and the condition \eqref{ineq:ellipticity} in Appendix~\ref{app:viscosity}} is obviously satisfied since the matrix $\sigma \sigma^\top$ is positive semidefinite.
	We first show that the value function $V^{\sf s}$ is a viscosity subsolution of \eqref{eq:SHJB}. 
	\mnr{For $p\ge 2$ satisfying \eqref{ineq:g_growth},} fix any 
	$\phi \in \mnr{C_p^{1,2}}([0,T]\times \bbR^n)$,
	and let $(\tx)$ be a \ito{global} maximum point of $V^{\sf s} - \phi$. 
\modi{Let us consider a constant control $u_s = \bar{u}$ for any $s\in [t,T]$, with $\bar{u}\in \bbU$.} Denote the corresponding state process $\state{x}_{s}^{t,x,\mnr{u}}$ by $\bar{\state{x}}_{s}$.
	Then, for $\tau\in (t,T)$, we have
	\begin{equation}\label{ineq:implicit}
	\bbE\left[ \phi(\tx) - \phi(\tau, \bar{\state{x}}_\tau) \right] \le \bbE\left[ V^{\sf s}(\tx) - V^{\sf s}(\tau, \bar{\state{x}}_\tau) \right] . 
	\end{equation}
	By using Lemma \ref{lem:DP_s}, we obtain
	\begin{align*}
	V^{\sf s}(\tx) &\le \bbE\left[ \int_t^\tau \psi_0 (\mnr{u}_s)ds + V^{\sf s}(\tau, \bar{\state{x}}_\tau) \right] \\
	&= (\tau - t) \psi_0 (\mnr{\bar{u}}) + \bbE\left[ V^{\sf s}(\tau, \bar{\state{x}}_\tau) \right] .
	\end{align*}
	Therefore,
	\[
	\bbE\left[ \phi(\tx) - \phi(\tau, \bar{\state{x}}_\tau) \right] \le (\tau -t) \psi_0 (\mnr{\bar{u}}) .
	\]
	Note that under the growth condition \eqref{ineq:phi_dphi_growth}, it holds that
	\begin{align*}
	&\lim_{\tau \searrow t} \frac{\bbE[\phi(\tau, \bar{\state{x}}_\tau)] - \phi(\tx)}{\tau - t}= D_x \phi(\tx) \cdot f(x,\mnr{\bar{u}}) \\
	&\qquad\qquad  + \frac{1}{2}\mathrm{tr} \left(\sigma \sigma^\top (x,\mnr{\bar{u}}) D_x^2 \phi(\tx)\right) +  \phi_t (\tx) 
	\end{align*}
	where It{\^o}'s formula is applied \cite{Fleming2006}.
	Therefore, we get
	\begin{align*}
	-D_x \phi(\tx) \cdot f(x,\mnr{\bar{u}}) -& \frac{1}{2}\mathrm{tr} \left(\sigma \sigma^\top (x,\mnr{\bar{u}}) D_x^2 \phi(\tx)\right) \\
	&\qquad - \phi_t (\tx) \le \psi_0 (\mnr{\bar{u}}) .
	\end{align*}
This inequality holds for all $\mnr{\bar{u}}\in \bbU$. This means
	\begin{align*}
	&-\phi_t(\tx) + H^{\sf s}(x,D_x\phi(\tx),D_x^2\phi(\tx))\le 0 .
	\end{align*}
	
	We next show that $V^{\sf s}$ is a viscosity supersolution of \eqref{eq:SHJB}. 
	\ito{Fix any $\phi \in \mnr{C_p^{1,2}}(\mnr{[0,T]}\times \bbR^n)$, and let $(\tx)$ be a \ito{global} minimum point of $V^{\sf s} - \phi$.
		Then,  for any $\varepsilon > 0$ and $\tau\in (t,T)$, by Lemma~\ref{lem:DP_s}, there exists $\tilde{u} \in \calU^{\sf s}[t,T]$, which depends on $\varepsilon$ and $\tau$, such that
		\modi{
		\begin{equation}
			V^{\sf s}(\tx) + (\tau - t)\varepsilon \ge \bbE\left[ \int_t^\tau \psi_0 (\tilde{u}_s) ds + V^{\sf s} (\tau, \tilde{\state{x}}_\tau) \right] ,
		\end{equation}
		where we denote $\state{x}_s^{t,x,\tilde{u}}$ by $\tilde{\state{x}}_s$. Therefore, it holds that
	}
		\begin{align}
		0&\ge \bbE \left[ V^{\sf s}(\tx) - \phi(\tx) - V^{\sf s} (\tau,\tilde{ \state{x}}_{\tau}) + \phi (\tau, \tilde{\state{x}}_{\tau}) \right] \nonumber\\
		&\ge -(\tau-t)\varepsilon + \bbE \left[ \int_t^{\tau} \psi_0 (\tilde{u}_s) ds +  \phi (\tau, \tilde{\state{x}}_{\tau})- \phi(\tx) \right] . \label{ineq:mini_ineq}
		\end{align}
		By applying It{\^o}'s formula, we obtain
		\begin{align}
		\bbE[\phi(\tx) - \phi(\tau, \tilde{\state{x}}_\tau) ] &= \bbE\Biggl[ -\int_t^\tau D_x \phi(s, \tilde{\state{x}}_s)\cdot f(\tilde{\state{x}}_s, \tilde{u}_s)ds \nonumber\\
		&  -\int_t^\tau \frac{1}{2} {\rm tr}(\sigma \sigma^\top \left(\tilde{\state{x}}_s, \tilde{u}_s) D_x^2 \phi(s,\tilde{\state{x}}_s)\right) ds \nonumber\\
		& -\int_t^\tau \phi_t (s, \tilde{\state{x}}_s) ds \Biggr] .\label{eq:dynkin}
		\end{align}
		Here, note that
		\begin{align}
		&\bbE\left[ \int_t^\tau D_x \phi(s,\tilde{\state{x}}_s)\cdot f(\tilde{\state{x}}_s, \tilde{u}_s) ds \right] \nonumber\\
		&\qquad = \bbE\left[ \int_t^\tau D_x \phi(\tx)\cdot f(x,\tilde{u}_s) ds\right] + o(\tau -t)  \label{eq:smallo_stoc} .
		\end{align}

		\modi{To see this, rewrite \eqref{eq:smallo_stoc} as
			\begin{align}
			&\underbrace{\bbE\biggl[\int_{t}^{\tau}\big\{D_x\phi(s, \tilde{\state{x}}_s)-D_x\phi(\tx)\big\} \cdot f(\tilde{\state{x}}_s,\tilde{u}_s)ds\biggr]}_{\triangleq I_1 (\tau)} \nonumber\\
			& \mspace{10mu} + \underbrace{\bbE\biggl[ \int_{t}^{\tau} D_x\phi(\tx) \cdot \big\{f(\tilde{\state{x}}_s,\tilde{u}_s)-f(x,\tilde{u}_s)\big\} ds \biggr]}_{\triangleq I_2(\tau)}  = o(\tau -t) . \label{eq:order_assume}
			\end{align}
		The first term $I_1$ is bounded above as follows.
		\begin{align*}
		&I_1(\tau) \le \bbE\left[ \int_t^\tau \|D_x \phi(s,\tilde{x}_s) - D_x \phi (t,x) \| \cdot \| f(\tilde{x}_s, \tilde{u}_s) \| ds \right] \\
		&\le \left\{ LK_1 (1+\|x\|) + K_f \right\}  \\
		&\qquad\qquad \times \bbE\left[ \int_t^\tau \|D_x \phi(s,\tilde{x}_s) - D_x \phi (t,x) \| ds \right] \\
		&\le \left\{LK_1(1+\|x\|) + K_f \right\} \\
		&\qquad\qquad \times  (\tau - t) \sup_{s\in [t,\tau]} \bbE\left[\|D_x \phi(s,\tilde{x}_s) - D_x \phi (t,x) \|  \right],
		\end{align*}
		where $L$ and $K_1$ satisfies \eqref{eq:Lip_fsigma} and \eqref{ineq:moment_upper}, respectively, and $K_f$ is some constant satisfying $\|f(0,u)\| \le K_f$ for all $u\in \bbU$. If it holds that
		\begin{equation}\label{eq:omit_part}
			\lim_{s \searrow t} \bbE\left[\|D_x \phi(s,\tilde{x}_s) - D_x \phi (t,x) \|  \right]=0  ,
		\end{equation}
		then we obtain $\lim_{\tau \searrow t} I_1 (\tau)/(\tau -t) = 0$.
		Indeed, we can show \eqref{eq:omit_part} under the condition $\phi \in C_p^{1,2} ([0,T]\times \bbR^n)$ along the same lines as the proof of Theorem 2 in \cite{Ito2020ifac}. Likewise, we get $\lim_{\tau \searrow t} I_2(\tau)/(\tau -t) = 0$ under Assumption~$(A_1)$, and therefore \eqref{eq:order_assume} holds.

		}

		By the same argument, we see that
		\begin{align*}
		&\bbE\left[\int_t^\tau \frac{1}{2} {\rm tr}(\sigma \sigma^\top \left(\tilde{\state{x}}_s, \tilde{u}_s) D_x^2 \phi(s,\tilde{\state{x}}_s)\right) ds\right] \\
		&\qquad =\bbE\left[\int_t^\tau \frac{1}{2} {\rm tr}(\sigma \sigma^\top \left(x, \tilde{u}_s) D_x^2 \phi(\tx)\right) ds\right] + o(\tau -t),
		\end{align*}
		\[
		\bbE\left[ \int_t^\tau \phi_t (s,\tilde{\state{x}}_s)ds \right] = (\tau -t)\phi_t (\tx) + o(\tau -t) .
		\]
		Then, it follows from \eqref{ineq:mini_ineq} and \eqref{eq:dynkin} that
		\begin{align*}
		-(\tau - t)\varepsilon &\le \bbE\Biggl[ \int_t^\tau \Bigl\{ -D_x \phi(\tx)\cdot f(x,\tilde{u}_s) \\
		&\quad - \frac{1}{2} {\rm tr}(\sigma \sigma^\top \left(x, \tilde{u}_s) D_x^2 \phi(\tx)\right) - \psi_0 (\tilde{u}_s) \Bigr\} ds \Biggr] \\
		&\quad -(\tau -t)\phi_t (\tx) +o(\tau - t) \\
		&\le (\tau -t) \sup_{u\in \bbU} \Bigl\{ -D_x \phi(\tx)\cdot f(x,u) \\
		&\qquad - \frac{1}{2} {\rm tr}(\sigma \sigma^\top \left(x, u) D_x^2 \phi(\tx)\right) - \psi_0 (u)\Bigr\} \\
		&\quad\qquad -(\tau -t)\phi_t (\tx) +o(\tau - t) .
		\end{align*}
		Divide both sides by $(\tau -t)$ and let $\tau \searrow t$, then
		\begin{align*}
		-\varepsilon &\le \sup_{u\in \bbU} \Bigl\{-D_x \phi(\tx)\cdot f(x,u) \\
		&\quad - \frac{1}{2} {\rm tr}(\sigma \sigma^\top \left(x, u) D_x^2 \phi(\tx)\right)  -\psi_0 (u) \Bigr\} - \phi_t (\tx) . 
		\end{align*}}The arbitrariness of $\varepsilon$ shows that $V^{\sf s}$ is a viscosity supersolution of \eqref{eq:SHJB}.
		\modi{Combining the above arguments with Lemma~\ref{lem:value_growth} completes the proof.}
	\hfill$\Box$
\end{pf}

\subsection{Optimality of a control}
Next, we provide a necessary condition and a sufficient condition for the $L^0$ optimality. The {\em second-order right parabolic superdifferential} \mnr{$D_{t+,x}^{1,2,+}$} and {\em subdifferential} \mnr{$D_{t+,x}^{1,2,-}$} are defined in Appendix~\ref{app:viscosity}.
\modi{The proof is same as the one of \cite[Chapter 5, Theorem 5.3, 5.7]{Yong1999}, noting that under assumptions $(A_1),(A_2)$, and $(A_3)$, Theorem~\ref{thm:conti_V_s}, \ref{thm:stochastic_viscosity}, and Lemma~\ref{lem:DP_s} hold.}
\modi{
\begin{lemma}\label{thm:opt_ctrl_stoc}
	Fix $T>0$ and $(\tx)\in [0,T) \times \bbR^n$.
	Assume $(A_1),(A_2)$, and $(A_3)$. \\
	{\bf (Necessary condition)}\\
	Let $(\Omega^*, \calF^*, \bbP^*, \{w_s^* \}, \{u_s^*\}) \in \calU^{\sf s}[t,T]$ be an optimal solution for Problem~\ref{prob:main_stochastic}, and $\{\state{x}_s^*\}$ be the corresponding optimal state trajectory. Then, for any
	\begin{align*}
	(q^*,p^*,M^*)&\in \calL_{\calF^*}^2 (t,T;\bbR^n) \times \calL_{\calF^*}^2 (t,T;\bbR^n) \\
	&\quad \times \calL_{\calF^*}^2 (t,T;\calS^n)
	\end{align*}
	satisfying
	\begin{equation}\label{eq:subdiff}
	(q_s^*, p_s^*, M_s^*) \in D_{t+,x}^{1,2,-} V^{\sf s}(s,\state{x}_s^*), \ {\rm a.e.} \ s\in [t,T], \ \bbP^* {\rm \mathchar`-a.s.},
	\end{equation}
	it must hold that
	\begin{equation}
	\bbE[ q_s^* ] \le \bbE\left[ G(\state{x}_s^*, u_s^*, p_s^*, M_s^*) \right]  ,\ {\rm a.e.} \ s\in [t,T] , \label{ineq:nec_cond}
	\end{equation}
	where we define
	\[
	G(x,u,p,M) \triangleq -f(x,u) \cdot p - \frac{1}{2} {\rm tr} \left(\sigma \sigma^\top (x,u) M\right) - \psi_0 (u) .
	\]
	{\bf (Sufficient condition)}\\
	Let $(\bar{\Omega}, \bar{\calF}, \bar{\bbP}, \{\bar{w}_s\}, \{\bar{u}_s\})\in \calU^{\sf s}[t,T]$ and $\{\bar{\state{x}}_s\}$ be the corresponding state trajectory. If there exists
	\[
	(\bar{q},\bar{p},\bar{M})\in \calL_{\bar{\calF}}^2 (t,T;\bbR^n) \times \calL_{\bar{\calF}}^2 (t,T;\bbR^n) \times \calL_{\bar{\calF}}^2 (t,T;\calS^n)
	\]
	satisfying
	\begin{equation}
	(\bar{q}_s, \bar{p}_s, \bar{M}_s) \in D_{t+,x}^{1,2,+}V^{\sf s}(s,\bar{\state{x}}_s), \ {\rm a.e.} \ s\in [t,T], \ \bar{\bbP}{\rm \mathchar`-a.s.},
	\end{equation}
	and
	\begin{align}\label{eq:suff_stoc}
	\bar{q}_s &= G(\bar{\state{x}}_s,\bar{u}_s, \bar{p}_s, \bar{M}_s) \nonumber\\
	&= \max_{u\in \bbU} G(\bar{\state{x}}_s,u,\bar{p}_s,\bar{M}_s), \ {\rm a.e.} \ s\in [t,T], \ \bar{\bbP}{\rm \mathchar`-a.s.},
	\end{align}
	then $\{\bar{u}_s\}$ is an optimal control process.
	\hfill$\lhd$
\end{lemma}
}

\modi{
Compared to the verification theorem~\cite[]{Fleming2006} that is well known as an optimality condition for the case when the value function is smooth, the above conditions are quite complicated and do not show explicitly the relationship between the optimal control value and the state value at the current time via the value function.
In view of this, we derive a novel necessary and sufficient condition that is similar to the verification theorem and therefore much clearer.
Now we introduce some assumptions:
\begin{itemize}
	\item[$(B_1)$] For any $u\in \calU^{\sf s}[t,T]$, the value function $V^{\sf s}$ defined by \eqref{eq:value_func}, admits $V_t^{\sf s}, D_x V^{\sf s}$, and $D_x^2 V^{\sf s}$ at $(s, \state{x}_s)$ almost everywhere $s\in [t,T]$ and almost surely;
	\item[$(B_2)$] For any $\rho \in \{V_t^{\sf s}, D_x V^{\sf s}, D_x^2 V^{\sf s}\}$, there exists a function $\varphi:[t,T] \rightarrow S \ (S=\bbR, \bbR^n, \calS^n)$ such that, for any $s\in [t,T]$,
	\begin{equation}
	\rho_{\varphi,s} (\state{x}) \triangleq \left\{
	\begin{array}{ll}
	\rho (s,\state{x}), & ~~ {\rm if}~ \rho ~{\rm exists~at} ~ (s,\state{x}) , \\
	\varphi(s), & ~~ {\rm otherwise} ,
	\end{array}
	\right. \ x\in \bbR^n 
	\end{equation}
	is Borel measurable;
	\item[$(B_3)$] For any $\rho \in \{V_t^{\sf s}, D_x V^{\sf s}, D_x^2 V^{\sf s}\}$, there exist constants $K>0$ and $p\ge 2$ such that
	\[ \|\rho (s,x)\| \le K(1+\|x\|^p) \]
	holds at any $(s,x)\in [t,T]\times \bbR^n$ where $\rho(s,x)$ exists.
\end{itemize}
The validity of the above 
assumptions is discussed in Remark~\ref{rmk:reasonable}.

\begin{theorem}\label{thm:feedback}
	Fix $T>0$ and $(\tx)\in [0,T)\times \bbR^n$. Assume $(A_1), (A_2), (A_3)$ and $(B_1), (B_2), (B_3)$. 
	Then, $(\Omega, \calF, \bbP, \{w_s\}, \{u_s\}) \in \calU^{\sf s}[t,T]$ is an optimal solution for Problem~\ref{prob:main_stochastic} if and only if
	\begin{align}\label{eq:optimal_iff}
		&u_s \in \argmax_{u\in\mathbb{U}} \left\{G\left(\state{x}_s,u,D_x V^{\sf s} (s,\state{x}_s), D_x^2 V^{\sf s} (s,\state{x}_s) \right) \right\} \nonumber\\
		&\hspace{4cm}\ {\rm a.e.} \ s\in [t,T], \ \bbP{\rm \mathchar`-a.s.},
	\end{align}
	where  $\{\state{x}_s\}$ is the corresponding state trajectory.
\end{theorem}
\begin{pf}
	For a given $(\Omega, \calF, \bbP, \{w_s\}, \{u_s\}) \in \calU^{\sf s}[t,T]$ and the corresponding state trajectory $\{x_s\}$, define a stochastic process
	\begin{equation}
		q_s \triangleq \left\{
		\begin{array}{ll}
			V_t^{\sf s} (s,\state{x}_s), & ~~ {\rm if}~ V_t^{\sf s} ~{\rm exists~at} ~ (s,\state{x}_s) \\
			\varphi(s), & ~~ {\rm otherwise} ,
		\end{array}
		\right.
	\end{equation}
	where $\varphi: [t,T]\rightarrow \bbR$ satisfies Assumption~$(B_2)$, which ensures that $\{q_s\}$ is an $\{\calF_s\}_{s\ge t}$-adapted process.
	By Assumption~$(B_1)$, it holds that $q_s = V_t^{\sf s}(s,\state{x}_s)$ almost everywhere $s\in [t,T]$ and $\bbP$-almost surely, and by a slight abuse of notation, we denote $q_s = V_t^{\sf s}(s,\state{x}_s)$. 
	In addition, Assumption~$(B_3)$ and Lemma~\ref{lem:estimate_moment} imply that $\bbE[\int_t^T \|q_s\|^2 ds] < +\infty$. To sum up, we get $\{q_s\} \in \calL_{\calF}^2(t,T; \bbR^n)$. Similarly, take $p_s = D_x V^{\sf s}(s,\state{x}_s)$ and $M_s = D_x^2 V^{\sf s}(s,\state{x}_s)$.
	Note that, by \eqref{eq:subdiff_property} in Appendix~\ref{app:viscosity}, it holds that
	\begin{align*}
		&(q_s, p_s, M_s) \in D_{t+,x}^{1,2,-} V^{\sf s}(s, \state{x}_s) \cap D_{t+,x}^{1,2,+} V^{\sf s}(s, \state{x}_s), \\ 
		&\hspace{4cm}{\rm a.e.}\ s\in [t,T], \ \bbP{\rm \mathchar`-a.s.}
	\end{align*}
	If $(\Omega, \calF, \bbP, \{w_s\},  \{u_s\})$ is an optimal solution, by Lemma~\ref{thm:opt_ctrl_stoc}, \eqref{ineq:nec_cond} holds for $q_s^* = q_s, p_s^* = p_s, M_s^* = M_s$, that is,
	\begin{align}
	&\bbE[ V^{\sf s}_t (s, \state{x}_s) ] \le \bbE\left[ G\left(\state{x}_s, u_s, D_x V^{\sf s} (s, \state{x}_s), D_x^2 V^{\sf s} (s, \state{x}_s )\right) \right]  ,\nonumber\\ 
	&\hspace{5cm}{\rm a.e.} \ s\in [t,T] . \label{ineq:Vt<G}
	\end{align}
	On the other hand, since $V^{\sf s}$ is a viscosity solution to the HJB equation \eqref{eq:SHJB} under assumptions $(A_1),(A_2)$, and $(A_3)$ by Theorem~\ref{thm:stochastic_viscosity}, $V^{\sf s}$ satisfies \eqref{eq:SHJB} at any point where $V_t^{\sf s}, D_x V^{\sf s}$, and $D_x^2 V^{\sf s}$ exist. Together with Assumption~$(B_1)$,
	\begin{align}
	\hspace{-0.3cm}V^{\sf s}_t (s, \state{x}_s) &= H^{\sf s}(\state{x}_s, D_x V^{\sf s}(s,\state{x}_s), D_x^2 V^{\sf s}(s,\state{x}_s)) \nonumber\\
	&= \sup_{u\in \bbU} G\left(\state{x}_s, u, D_x V^{\sf s} (s,\state{x}_s), D_x^2 V^{\sf s} (s,\state{x}_s)\right) \nonumber\\
	&\ge G\left(\state{x}_s, u_s, D_x V^{\sf s} (s,\state{x}_s), D_x^2 V^{\sf s} (s,\state{x}_s)\right), \nonumber\\
	&\hspace{3cm} {\rm a.e.} \ s\in [t,T], \ \bbP{\rm \mathchar`-a.s.}  \label{ineq:Vt>G}
	\end{align}
	Combining \eqref{ineq:Vt<G} and \eqref{ineq:Vt>G}, we obtain
	\begin{align}\label{eq:nec_stoc}
	\hspace{-0.3cm}V^{\sf s}_t (s,\state{x}_s) &= \max_{u\in \bbU} G\left(\state{x}_s, u, D_x V^{\sf s} (s,\state{x}_s), D_x^2 V^{\sf s} (s,\state{x}_s)\right) \nonumber\\
	&= G\left(\state{x}_s, u_s, D_x V^{\sf s} (s,\state{x}_s), D_x^2 V^{\sf s} (s,\state{x}_s)\right), \nonumber\\
	&\hspace{3cm} {\rm a.e.} \ s\in [t,T], \ \bbP{\rm \mathchar`-a.s.} 
	\end{align}
	Hence, the optimal control process $\{u_s\}$ satisfies \eqref{eq:optimal_iff}.
	
	Conversely, if an admissible control process $\{u_s\}$ satisfies \eqref{eq:optimal_iff}, it follows from \eqref{ineq:Vt>G} that
	\begin{align*}
	q_s &= \sup_{u\in \bbU} G(\state{x}_s,u,p_s,M_s) \nonumber\\
	&= G(\state{x}_s,u_s, p_s, M_s), \ {\rm a.e.} \ s\in [t,T], \ \bbP{\rm \mathchar`-a.s.}  \nonumber\\
	\end{align*}
	 In other words, the sufficient condition \eqref{eq:suff_stoc} in Lemma~\ref{thm:opt_ctrl_stoc} holds for $\bar{q}_s = q_s, \bar{p}_s = p_s, \bar{M}_s = M_s, \bar{\state{x}}_s = \state{x}_s, \bar{u}_s = u_s$. This completes the proof.
	 \hfill$\Box$
\end{pf}

Theorem~\ref{thm:feedback} can be seen as a generalization of the verification theorem for almost everywhere differentiable value functions; see Remark~\ref{rmk:reasonable} below.
It should be emphasized that the above result can be easily generalized to the cost functional with a uniformly continuous state cost and control cost since Theorem 5.3 and 5.7 in \cite[Chapter 5]{Yong1999} hold.

\begin{rmk}\label{rmk:reasonable}
	\modi{If $V^{\sf s}$ admits any $\rho \in \{V_t^{\sf s}, D_x V^{\sf s}, D_x^2 V^{\sf s}\}$ almost everywhere $(s,x) \in [t,T]\times \bbR^n$, and $\rho$ is continuous almost everywhere, then, from Lusin's theorem, there exists a Borel measurable function which coincides with $\rho$ almost everywhere. Hence, in this case, we can remove Assumption~$(B_2)$. In addition, if for any $u\in \calU^{\sf s}[t,T]$, there exist densities of $\{x_s\}$, Assumption~$(B_1)$ holds. This is a sufficient condition, but it is not necessary. See \cite[]{Bouleau2010} for the existence of the density for solutions to stochastic differential equations. For Assumption~$(B_3)$, we expect that the condition can be removed or relaxed along the line of \cite[Theorem IV.3.1]{Fleming2006}, but this is not our focus in the present paper.
	\hfill$\lhd$}
\end{rmk}

	Theorem~\ref{thm:feedback} immediately characterizes the $L^0$ optimal control in terms of the feedback control. 
	In fact, as a straightforward consequence of Theorem~\ref{thm:feedback}, we obtain the following result.
	\begin{corollary}\label{cor:state-feedback}
	Fix $T>0$ and $(t,x)\in [0,T)\times \bbR^n$.
	Assume $(A_1), (A_2), (A_3)$ and $(B_1), (B_2), (B_3)$.
	Let a Borel measurable function $\underline{u}:[t,T]\times \bbR^n \rightarrow \bbU$ satisfy
	\begin{align*}
		&\underline{u}(s,x^{\prime}) \in \argmax_{u\in\mathbb{U}} \left\{G\left(\state{x}^{\prime},u,D_x V^{\sf s} (s,\state{x}^{\prime}), D_x^2 V^{\sf s} (s,\state{x}^{\prime}) \right) \right\} \\
		&\hspace{5cm} {\rm a.e.} \ s\in [t,T] , \ x^{\prime} \in \bbR^n.
	\end{align*}
	Fix any reference probability space $(\Omega, \calF, \bbP, \{w_s\})$. If the stochastic differential equation
	\begin{align*}
		&dx_s = f(x_s, \underline{u}(s,x_s)) ds + \sigma(x_s, \underline{u}(s,x_s)) dw_s, \ s>t, \\
		&x_t = x
	\end{align*}
	has a unique strong solution, then $u_s^* \triangleq \underline{u}(s,x_s)$ is an optimal control process, namely, $(\Omega, \calF, \bbP, \{w_s\}, \{u_s^*\}) \in \calU^{\sf s}[t,T]$ is an optimal solution for Problem~\ref{prob:main_stochastic}.
	\hfill$\lhd$
	\end{corollary}
	Here, we emphasize that in the above result, we can choose any reference probability space to be fixed. Thus, for a state-feedback controller, we need not to distinguish which reference probability space is optimal, and we can concentrate only on control processes.

}

\section{Characterization of \rev{sparse optimal stochastic} control}\label{sec:char_stochastic}
In this section, we focus on the control-affine systems satisfying
\begin{equation}\label{eq:affine}
f(\state{x},u) = f_0 (\state{x}) + \sum_{j=1}^m f_j(\state{x}) u^{(j)} , \ \sigma(\state{x},u) = \sigma(\state{x})
\end{equation}
for some $f_j : \bbR^n \rightarrow \bbR^n, j = 0,1,2,\ldots, m$ where $u^{(j)}$ is the $j$-th component of $u\in \bbR^m$.
First, we reveal the discreteness of the stochastic $L^0$ optimal control.
Next, we show an equivalence between the $L^0$ optimality and the $L^1$ optimality. Thanks to the equivalence, we \mnr{ensure} that our value function is a classical solution of the associated HJB equation under some assumptions.

\subsection{Discreteness of the optimal control}
We explain the discreteness of the stochastic $L^0$ optimal control based on \mnr{Theorem~\ref{thm:feedback}}.

\mnr{
\begin{theorem}\label{thm:discrete}
	Fix $T>0$ and $(t,x)\in [0,T)\times \bbR^n$. Assume $(A_1),(A_2),(A_3)$ and \mnr{$(B_1), (B_2), (B_3)$}. If the system \eqref{eq:stochastic_system} is control-affine, i.e., \eqref{eq:affine} holds, and $\mathbb{U}=\{u\in\mathbb{R}^m: U_{j}^{-} \leq u^{(j)} \leq U_{j}^{+}, \forall j\}$ for some $U_{j}^{-}< 0$ and $U_{j}^{+}>0$, then, $u \in \calU^{\sf s}[t,T]$ is an optimal solution to Problem~\ref{prob:main_stochastic} if and only if
	\begin{align}
	u_s^{(j)} &\in \argmax_{u^{(j)}\in \bbU_j} \{ -(f_j (\state{x}_s) \cdot D_x V^{\sf s}(s,\state{x}_s))u^{(j)} - |u^{(j)}|^0 \} \nonumber\\
	&\hspace{4cm} {\rm a.e.} \ s\in [t,T], \ \bbP{\rm \mathchar`-a.s.}  \label{eq:iff_affine}
	\end{align}
	for all $j=1,2,\ldots,m$ where $\bbU_j \triangleq \{a\in \bbR : U_j^- \le a \le U_j^+\}$.
	Furthermore, if an optimal control process $\{u_s^*\}$ exists, then the $j$-th component of $u_s^*$ takes only three values of $\{U_j^-, 0, U_j^+\}$ almost everywhere $s\in [t,T]$ and almost surely.
\end{theorem}
\begin{pf}
	By \mnr{Theorem~\ref{thm:feedback}}, a necessary and sufficient condition for the $L^0$ optimality of $\{u_s\}$ is given by
	\begin{align}
	u_s &\in \argmax_{u\in\mathbb{U}} G(\state{x}_s, u, D_x V^{\sf s}(s,\state{x}_s), D_x^2 V^{\sf s}(s,\state{x}_s)) \nonumber\\
	&= \argmax_{u\in\mathbb{U}} \left\{- \sum_{j=1}^m f_j (\state{x}_s)u^{(j)} \cdot  D_xV^{\sf s} (s,\state{x}_s) \right. \nonumber\\
	&\qquad\qquad\qquad \left. - \sum_{j=1}^m |u^{(j)}|^0 \right\} \nonumber\\
	&= \argmax_{u\in\mathbb{U}} \sum_{j=1}^m \left\{-  (f_j (\state{x}_s)\cdot D_xV^{\sf s} (s,\state{x}_s)) u^{(j)} \right. \nonumber\\
	&\qquad\qquad\qquad \left. -  |u^{(j)}|^0 \right\}, \ {\rm a.e.} \ s\in [t,T], \ \bbP{\rm \mathchar`-a.s.} \label{eq:iff_affine1}
	\end{align}
	noting that $\sigma$ does not depend on the control variable.
	Then, \eqref{eq:iff_affine1} is equivalent to \eqref{eq:iff_affine}.
	
	Next, it follows from \eqref{eq:iff_affine} and an elementary calculation that
	\begin{equation}\label{eq:discrete_L0}
	u_s^{(j)}
	\in\begin{cases}
	\{U_{j}^{-}\} , &\mbox{if}\quad b_{j}(s,\state{x}_s)U_{j}^{-}<-1,\\
	\{U_{j}^{-}, 0\} , &\mbox{if}\quad  b_{j}(s,\state{x}_s)U_{j}^{-}=-1,\\
	\{0\} , &\mbox{if}\quad  b_{j}(s,\state{x}_s)U_{j}^{-}>-1,\\ 
	&\mspace{12mu}\quad\mbox{and~} b_{j}(s,\state{x}_s)U_{j}^{+}>-1,\\
	\{0,U_{j}^{+}\} , &\mbox{if}\quad  b_{j}(s,\state{x}_s)U_{j}^{+}=-1,\\
	\{U_{j}^{+}\} , &\mbox{if}\quad  b_{j}(s,\state{x}_s)U_{j}^{+}<-1,
	\end{cases}
	\end{equation}
	where we define $b_j (s,\state{x}) \triangleq D_x V^{\sf s} (s,\state{x}) \cdot f_j (\state{x})$. Therefore, the $j$-th component of an optimal control process $\{u_s^*\}$ must take only three values of $\{U_j^-, 0, U_j^+\}$ almost everywhere $s\in [t,T]$ and almost surely.
	\hfill$\Box$
\end{pf}

}

\subsection{Equivalence between $L^0$ optimality and $L^1$ optimality}
\modi{Let us consider the stochastic $L^1$ optimal control problem where the cost functional $J^{\sf s}$ in Problem~\ref{prob:main_stochastic} is replaced by the following one:
\begin{equation}\label{eq:L1_cost}
	J_1^{\sf s} (t,x,u) \triangleq \bbE \left[\sum_{j=1}^m \int_t^T |u_{s}^{(j)}| ds + g(\state{x}_T^{t,x,u}) \right] .
\end{equation}
The corresponding value function is defined by
\begin{equation}\label{eq:L1_value_func}
V_1^{\sf s} (\tx)\triangleq \inf_{u\in \calU^{\sf s}[t,T]} J_1^{\sf s} (t,x,u) ,\ (t,x) \in [0,T]\times \bbR^n .
\end{equation}

We here show the coincidence of the value functions of the $L^0$ optimal control and the $L^1$ optimal control for the control-affine system.

\begin{theorem}\label{cor:L0_L1_s}
	Fix $T>0$ and $(t,x)\in [0,T)\times \bbR^n$. Assume $(A_1), (A_2)$, and $(A_3)$. If the system \eqref{eq:stochastic_system} is control-affine, i.e., \eqref{eq:affine} holds, and $\bbU = \{u\in \bbR^m : |u^{(j)}| \le 1, \forall j\}$, then for the value functions $V^{\sf s}$ and $V_1^{\sf s}$ defined by \eqref{eq:value_func} and \eqref{eq:L1_value_func}, respectively, it holds that
	\[
	 V^{\sf s}(\tx) = V_1^{\sf s}(\tx), \  \forall (\tx) \in [0,T]\times \bbR^n.
	\]
	\modi{In addition, $V^{\sf s}$ is a unique, at most polynomially growing viscosity solution to the HJB equation \eqref{eq:SHJB} with the terminal condition \eqref{eq:SHJB_initial}.}
\end{theorem}  }
\begin{pf}
	In this setting, for any $x,p\in \bbR^n$ and $M\in \calS^n$,
	\begin{align*}
	H^{\sf s}(x,p,M) &= \sup_{u\in \bbU} \left\{ -\sum_{j=1}^m f_j (x)u^{(j)} \cdot p - \sum_{j=1}^m |u^{(j)}|^0 \right\}\\
	&\qquad  - f_0 (x) \cdot p - \frac{1}{2} {\rm tr}(\sigma \sigma^\top (x)M) \\
	&= \sum_{j=1}^m \sup_{u^{(j)} \in \bbU_j} \left\{ -(f_j (x) \cdot p) u^{(j)} - |u^{(j)}|^0 \right\} \\
	&\qquad - f_0 (x) \cdot p - \frac{1}{2} {\rm tr}(\sigma \sigma^\top (x)M)
	\end{align*}
	where $\bbU_j = \{a\in \bbR : |a| \le 1\}$. 
	Here, it follows from an elementary calculation that
	\begin{align*}
	&\underset{u^{(j)} \in\mathbb{U}_j}{\sup}\bigg\{- a_{x,p}^{j} u^{(j)} - |u^{(j)}|^0\bigg\} \\
	&\qquad = \underset{u^{(j)}\in\mathbb{U}_j}{\sup}\bigg\{- a_{x,p}^{j} u^{(j)} - |u^{(j)}|\bigg\}
	\end{align*}
	for all $x,p\in\mathbb{R}^n$ and $j=1,2,\dots,m$,
	where $a_{x,p}^{j}\triangleq f_j(x)\cdot p$.
	Indeed, the supremum of both sides is given by 
	\begin{align*}
	\begin{cases}
	a_{x,p}^{j} - 1 , &\mbox{if}\quad a_{x,p}^{j}>1,\\
	0 , &\mbox{if}\quad |a_{x,p}^{j}| \le 1,\\ 
	-a_{x,p}^{j} - 1 , &\mbox{if}\quad a_{x,p}^{j}<-1.\\
	\end{cases}
	\end{align*}
	Hence, the HJB equation \eqref{eq:SHJB} is equivalent to
	\begin{equation}
	-v_t(\tx) + H_1(x,D_x v(\tx), D_x^2 v(\tx)) = 0,
	\label{eq:SHJB_L1}
	\end{equation}
	where 
	\begin{align*}
	H_1(x,p,M) &\triangleq\sup_{u\in\mathbb{U}} \Bigl\{-f(x,u) \cdot p -\frac{1}{2} {\rm tr}(\sigma \sigma^\top (x,u) M)  \\
	&\qquad - \psi_1(u)\Bigr\}, \quad  x, p\in\mathbb{R}^n, M\in \calS^n , \\
	\psi_1(a) &\triangleq\sum_{j=1}^{m} |a_j|, \quad a\in\mathbb{R}^m.
	\end{align*}
	Note that the equation \eqref{eq:SHJB_L1} is the HJB equation for the $L^1$ optimal control problem.
	Moreover, it is known that \modi{$V_1^{\sf s}$ defined via the $L^1$ optimal control} is a unique, at most polynomially growing viscosity solution to the associated HJB equation \eqref{eq:SHJB_L1} \modi{or, equivalently \eqref{eq:SHJB}} with the terminal condition \eqref{eq:SHJB_initial} \cite[]{Pham2009}.

	Now, $V^{\sf s}$ is also a viscosity solution to the HJB equation \eqref{eq:SHJB} with the terminal condition \eqref{eq:SHJB_initial} by Theorem~\ref{thm:stochastic_viscosity}.
	\ito{Note also that $V^{\sf s}$ satisfies the polynomial growth condition by Lemma~\ref{lem:value_growth}.} \modi{Therefore, by the aforementioned uniqueness of the viscosity solution to the HJB equation \eqref{eq:SHJB}, we conclude that $V^{\sf s} = V_1^{\sf s}$.}
	\hfill$\Box$
\end{pf}

\modi{Theorem~\ref{cor:L0_L1_s} justifies the use of the value function for the $L^1$ optimal control to obtain the $L^0$ optimal control. For example, we can use a sampling-based algorithm recently proposed in \cite{ExaThe18} to calculate the value function.
}

In contrast to the deterministic case where the corresponding HJB equation is of first order, if the second order HJB equation is uniformly elliptic, then we expect that the HJB equation with a terminal condition has a unique classical solution. By using this property and Theorem~\ref{cor:L0_L1_s}, we \mnr{show} that the value function $V^{\sf s}$ is a unique classical solution to the HJB equation under some assumptions. Define
\begin{align*}
C_b^k(\bbR^n) \triangleq \{ \rho \in C^k(\bbR^n) : \rho \textrm{~and~all~partial~derivatives~of~}\\
\rho \textrm{~of~orders}\le k\textrm{~are~bounded} \} .
\end{align*}
\begin{corollary}
	Suppose the assumptions in Theorem~\ref{cor:L0_L1_s} and the following assumptions:
	\begin{itemize}
		\item[$(a)$] For any $\rho \in\{f_0, f_1 \ldots f_m,\sigma\sigma^\top\}$, $\rho \in C_b^2(\bbR^n)$,
		\item[$(b)$] $g\in C_b^3(\bbR^n)$,
		\item[$(c)$] {\bf Uniform ellipticity condition:}\\
		There exists $c>0$ such that, for all $x\in \bbR^n$ and $\xi \in \bbR^n$,
		\[ \xi^\top \sigma\sigma^\top (x) \xi \ge c\|\xi\|^2. \]
	\end{itemize}
	Then, the value function $V^{\sf s}$ is a unique classical solution to the HJB equation \eqref{eq:SHJB} with the terminal condition \eqref{eq:SHJB_initial}.
\end{corollary}
\begin{pf}
	By \cite[Theorem IV.4.2]{Fleming2006}, the HJB equation \eqref{eq:SHJB_L1} with the terminal condition \eqref{eq:SHJB_initial} for the $L^1$ optimal control problem has a bounded unique classical solution \modi{under assumptions~$(a), (b), (c)$}.
	In other words,
	the HJB equation \eqref{eq:SHJB} with \eqref{eq:SHJB_initial} has a bounded unique classical solution. \ito{Note that any classical solution of \eqref{eq:SHJB} is also a viscosity solution.}
	Note also that the value function $V^{\sf s}$ is a unique viscosity solution satisfying a polynomial growth condition by Theorem~\ref{cor:L0_L1_s}. This means that $V^{\sf s}$ \mnr{must be} a unique classical solution to \eqref{eq:SHJB} with \eqref{eq:SHJB_initial}.
	\hfill$\Box$
\end{pf}

\modi{
Thanks to the above result, we need not to consider the non-differentiability of the value function, and we can apply usual numerical methods to solve the HJB equation under the conditions~$(a), (b), (c)$.
}

\modi{
	In Theorem~\ref{cor:L0_L1_s}, we have shown the equivalence about the value functions of the $L^0$ optimal control and the $L^1$ optimal control.
	Combining this and the discreteness of the $L^0$ optimal control, we obtain an equivalence for the optimal control itself.

\begin{corollary}
	Suppose the assumptions in Theorem~\ref{thm:discrete} and let $U_j^- = -1, U_j^+ = 1, j = 1,\ldots, m$. If an $L^0$ optimal control process exists, then it is also an $L^1$ optimal control process. Conversely, if an $L^1$ optimal control process $\{u_s^{1*}\}$ exists, and it holds that
	\begin{equation}\label{eq:normality}
		| f_j(\state{x}_s)\cdot D_x V^{\sf s}(s,\state{x}_s) | \neq 1, \ {\rm a.e.} \ s\in [t,T], \ \bbP{\rm \mathchar`-a.s.},
	\end{equation}
	where $\{\state{x}_s\}$ is the corresponding optimal state trajectory, then $\{u_s^{1*}\}$ is also an $L^0$ optimal control process.
\end{corollary}
\begin{pf}
	By Theorem~\ref{thm:discrete}, each element of an $L^0$ optimal control $\{u_s^*\}$ takes only three values of $\{-1,0,1\}$, and therefore it holds that $|u_s^{*^{(j)}}| = |u_s^{*^{(j)}}|^0$. In addition, the optimal values of \eqref{eq:cost_stoc} and \eqref{eq:L1_cost} coincide by Theorem~\ref{cor:L0_L1_s}. This implies that $\{u_s^*\}$ is an $L^1$ optimal control process. Next, by the same arguments as in the proofs of Theorem~\ref{thm:feedback} and \ref{thm:discrete}, a control process $\{u_s\}$ is an $L^1$ optimal control process if and only if
	\begin{align*}
	u_s^{(j)}
	&\in\begin{cases}
	\{-1\} , &\mbox{if}\quad b_{j}(s,\state{x}_s) > 1,\\
	[-1, 0] , &\mbox{if}\quad  b_{j}(s,\state{x}_s)=1,\\
	\{0\} , &\mbox{if}\quad   |b_{j}(s,\state{x}_s) | < 1,\\
	[0,1] , &\mbox{if}\quad  b_{j}(s,\state{x}_s)=-1,\\
	\{1\} , &\mbox{if}\quad  b_{j}(s,\state{x}_s)<-1,
	\end{cases}
	\ s\in [t,T], \ \bbP{\rm \mathchar`-a.s.},
	\end{align*}
	where $b_j (s,\state{x}) = f_j (\state{x}) \cdot  D_xV^{\sf s}(s,\state{x}) =  f_j (\state{x}) \cdot  D_xV_1^{\sf s}(s,\state{x})$. Therefore, if \eqref{eq:normality} holds, then each element of the $L^1$ optimal control also takes only three values of $\{-1,0,1\}$. Then we obtain the desired result as in the first part of the proof.
	\hfill$\Box$
\end{pf}
The condition~\eqref{eq:normality} corresponds to the {\em normality} of the $L^1$ optimal control problem~\cite{NagQueNes16}.

}

\modi{
\begin{rmk}
	Finally, we would like to point out that all the results obtained in this paper can be extended to the case where a continuous state transition cost $\ell:\bbR^n \rightarrow \bbR$ is added to our cost functional \eqref{eq:cost_stoc}, i.e.,
	\[
		J_{\ell}^{\sf s} (t,x,u) \triangleq \bbE\left[ \int_t^T \left( \ell(x_s) + \psi_0 (u_s) \right) ds + g(x_T) \right] .
	\]
	Indeed, since $\ell$ does not depend on $u$, the only difference in the associated HJB equation is an additional term $\ell (x)$. Moreover, the continuity of $\ell$ can be used to prove the continuity of the corresponding value function.
	\hfill$\lhd$
\end{rmk}
}

\begin{example}[Revisited] \label{ex:2}
\ikeda{
	\modi{Throughout the following examples, we fix a reference probability space and consider a state-feedback controller.}
	We explain the result for (\ref{eq:simple_s}) in more detail. 
	\ito{
		First, we consider the deterministic case, i.e., $\sigma = 0$. 
		We can show that a \rev{smooth} function having the polynomial growth property
		satisfies the associated HJB equation. By the uniqueness of the viscosity solution, this is the value function; \mnr{see Theorem}~\ref{cor:L0_L1_s}.
		Note also that it is possible to apply \cite[Theorem 4]{IkeKas_ECC19} without the Lipschitz continuity of $g$ due to the smoothness of $V^{\sf s}$.
		Therefore, it can be verified that the $L^0$ optimal feedback control $u^* (s,\state{x})$ is given by
		\begin{equation}\label{eq:determinic_law}
		u^{\ast}(s,\state{x}) = 
		\begin{cases}
		-1,  &\mbox{if~} \frac{1}{2}e^{-2c(T-s)}<\state{x}, ~0\leq s\leq T, \\
		0,  &\mbox{if~} |\state{x}|\le \frac{1}{2}e^{-2c(T-s)}, ~0\leq s\leq T, \\
		1,  &\mbox{if~} \state{x}<-\frac{1}{2}e^{-2c(T-s)}, ~0\leq s\leq T.
		\end{cases}
		\end{equation}

		This analysis implies that the Lipschitz continuity of $g$ is not necessary for the value function to be differentiable almost everywhere; see Remark~\ref{rmk:conti_Vs}.
		
		Next, we consider the stochastic case, i.e., $\sigma > 0$. 
		The associated HJB equation is
		\[
		\mspace{10mu}
		\begin{cases}
		-v_t(\tx) -cx D_x v(\tx) - \frac{\sigma^2}{2} D_x^2 v (\tx) \\
		\hspace{2cm} + \alpha(D_x v(\tx)) = 0, ~~ (\tx) \in [0,T)\times \mathbb{R}, \\
		v(T,x) = x^2, \quad  x\in \mathbb{R} ,
		\end{cases}
		\]
		\modi{where 
			\[
			\alpha (p) \triangleq
			\begin{cases}
			p-1,  &\mbox{if~} p\geq 1, \\
			0,  &\mbox{if~} |p|<1, \\
			-p-1,  &\mbox{if~} p\leq -1.\\
			\end{cases}
			\]}\rev{We solve the above HJB equation numerically using a finite difference scheme. See \cite[]{Fleming2006,ExaThe18} for numerical methods to compute the viscosity solution to the HJB equation.}
		We take $c = 1, \sigma = 0.1$, and $T = 1$. The switching \rev{boundary} $\{(s,\state{x}) : |D_x V^{\sf s} (s,\state{x})| = 1 \}$ is depicted in Fig.~\ref{fig:switching_curve}. For comparison, we also plot the deterministic optimal switching \rev{boundary} $\{ (s,\state{x}) : |\state{x}| =\frac{1}{2} e^{-2c(T-s)} \}$ obtained in (\ref{eq:determinic_law}). As shown in Fig.~\ref{fig:switching_curve}, the region where the stochastic $L^0$ optimal control takes value $0$ is larger than the deterministic one.
		This implies that the stochastic $L^0$ optimal control gives a \mnr{sparser} control than the deterministic one instead of allowing the larger variance of the terminal state. 
		\begin{figure}[tb]
			\centering
			\includegraphics[scale=0.4]{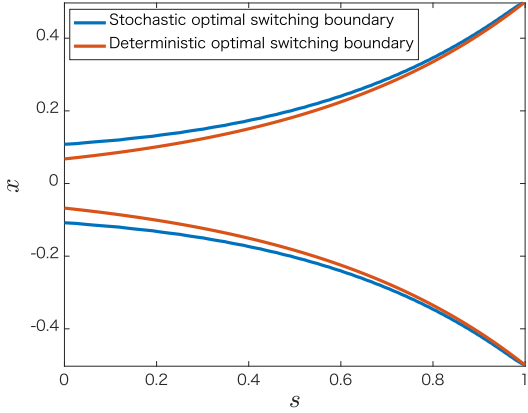}
			\caption{Stochastic optimal switching \rev{boundary} obtained by the numerical solution $V^{\sf s}$ (blue) and the deterministic optimal switching \rev{boundary} (\eqref{eq:determinic_law}, red).}
			\label{fig:switching_curve}
		\end{figure}
	}
}
\hfill$\lhd$
\end{example}

\modi{
\begin{example} \label{ex:3}
	Next, we consider a simplified load frequency control (LFC) model depicted in Fig.~\ref{fig:block_LFC}; see \cite{Kashima2014,Kashima2019} for more details.
	The physical meanings of $ x_s^{(1)} $ and $ x_s^{(2)} $ are frequency deviation and its compensation by a thermal plant, respectively. The feedback loop with $ 1/s $ and the saturation function
	\begin{equation}
	{\rm sat}_d (x) \triangleq
	\begin{cases}
	- d, &x < -d, \\
	x, & |x| \le d, \\
	d, &x > d,
	\end{cases}
	\hspace{0.5cm} x \in \bbR ,
	\end{equation}
	represents the rate limiter, where $d > 0$ characterizes the limited responsiveness of the adjustment of the thermal power generation.
	An extra compensation, which should not be activated for long time, is denoted as $ u_s $.
	The dynamics in Fig.~\ref{fig:block_LFC} is given by
	\begin{equation}\label{eq:nonlinear_ex}
		\begin{cases}
		d\state{x}_s^{(1)} &= (-p \state{x}_s^{(1)} - k\state{x}_s^{(2)} ) ds + ku_s ds + k\sigma dw_s, \\
		d\state{x}_s^{(2)} &= {\rm sat}_d (\state{x}_s^{(1)} - \state{x}_s^{(2)}) ds ,
		\end{cases}
	\end{equation}
	where $p > 0, k > 0, \sigma > 0$.
	We take $p = 1/3, k = 2, \sigma = 0.5, \bbU = [-1,1], T = 0.5$, and $g(x) = \|x\|^2$. Based on the equivalence result in Theorem~\ref{cor:L0_L1_s}, we employ a sampling-based method proposed in \cite[]{ExaThe18} with radial basis functions to solve the associated HJB equation. Figure~\ref{fig:boundary_initial} compares the obtained switching boundaries at time $s=0$, i.e., $\{x : k|(D_x V^{\sf s})^{(1)} (0,\state{x}) | = 1 \}$ for $d = 0.4$ and the linear case $(d = +\infty)$. To describe the result, let us consider the case $\state{x}^{(1)} > 0$ and $\state{x}^{(2)} \simeq 0$.
	In such a case, it is expected that $ x^{(2)} $ increases to suppress $ x^{(1)} $. When the rate limiter prevents the quick adjustment of $ x^{(2)} $, we need to activate $ u_s $. This is why the region on which the optimal control takes value $ 0 $ is larger for $ d=+\infty $ than for $ d=0.4 $.
	Similar interpretation applies to the case with $ x^{(1)} \simeq 0 $ and $ x^{(2)} > 0$.
	\hfill$\lhd$

	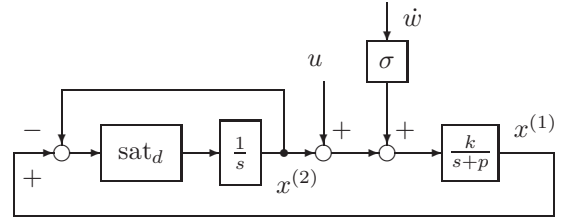
\begin{figure}[t]
		\setlength\unitlength{3pt}
		\begin{center}
			\begin{picture}(110,30)(10,0)

			\put(50,7){{$x^{(2)}$}}
			
			\put(17,4){\line(1,0){68}}
			\put(17,4){\line(0,1){8}}
			\put(85,4){\line(0,1){8}}

			\put(17,12){\vector(1,0){5}}
			\put(23,12){\circle{2}}
			\put(24,12){\vector(1,0){4}}
			\put(28,8.5){\framebox(10,7){${\rm sat}_{d}$}}
			\put(38,12){\vector(1,0){5}}
			\put(43,8){\framebox(5,8){$\frac{1}{s}$}}
			\put(48,12){\vector(1,0){7}}
			\put(51,12){\line(0,1){8}}
			\put(51,20){\line(-1,0){28}}
			\put(23,20){\vector(0,-2){7}}
			
			\put(56,12){\circle{2}}
			
			\put(57,12){\vector(1,0){6}}
			\put(64,12){\circle{2}}
			
			\put(56,21){\vector(0,-1){8}}
			\put(65,12){\vector(1,0){6}}

			\put(64,21){\vector(0,-1){8}}
			\put(61.5,21){\framebox(5,5){$\sigma$}}
			\put(64,31){\vector(0,-1){5}}
			\put(78,12){\line(1,0){7}}
			
			\put(71,8.5){\framebox(7,7){$\frac{k}{s+p}$}}

			\put(51,12){\circle*{1}}
			
			\put(66,28){$\dot w$}
			\put(54,23){$u$}
			
			\put(18,8){$+$}
			\put(18,14){$-$}
			
			\put(57,14){$+$}
			\put(65,14){$+$}
			
			\put(80,14){$x^{(1)}$}

			\end{picture}
		\end{center}
		\caption{Block diagram of the load frequency control system.}
		\label{fig:block_LFC}
	\end{figure}

	\begin{figure}[tb]
		\centering
		\includegraphics[scale=0.4]{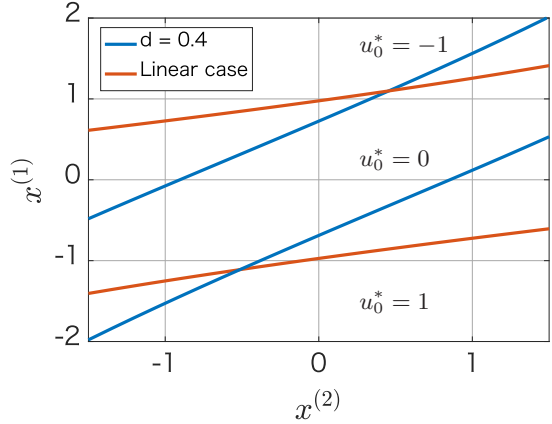}
		\put(-93,46){\small $u_0^* = 1$}
		\put(-93,143){\small $u_0^* = -1$}
		\put(-93,100){\small $u_0^* = 0$}
		\caption{Optimal switching boundaries at time $s = 0$ for $d = 0.4$ (blue) and the linear case (red), and the optimal control value $u_0^*$.}
		\label{fig:boundary_initial}
	\end{figure}

\end{example}
}

\section{Conclusions}
\label{sec:conclusion}
We have investigated a finite horizon stochastic optimal control problem with the $L^0$ control cost functional.
We have characterized the value function as a viscosity solution to the associated HJB equation 
and shown an equivalence theorem between the $L^0$ optimality and the $L^1$ optimality via \mnr{the uniqueness of a viscosity solution}. Thanks to the equivalence, we \mnr{have ensured that the value function is a classical solution of the associated HJB equation under some conditions}.
Moreover, we have derived a sufficient and necessary condition for the $L^0$ optimality 
that connects the current state and the current optimal control value. 
Furthermore, we have revealed the discreteness property of the sparse optimal stochastic control for control-affine systems.

\begin{ack}
	\rev{This work was supported in part by JSPS KAKENHI under Grant Number JP18H01461 and by JST, ACT-X under Grant Number JPMJAX1902.}
\end{ack}

\appendix
\vspace{\baselineskip}
{\bf Appendix}
\ike{\section{Moment estimate for the state}\label{app:moment}
	Here, we introduce an estimate for the $p$-th order moment of the state governed by the stochastic system \eqref{eq:stochastic_system} \cite[Theorem 1.2]{Nisio2014}.
	\begin{lemma}\label{lem:estimate_moment}
		Fix $T>0$.
		Assume $(A_1)$ and let $p\ge 2$ be given. Then there exists a positive constant $K_p$ such that, for any $(t,x)\in [0,T] \times \bbR^n$,
		\begin{equation}\label{ineq:moment_upper}
		\bbE\left[ \sup_{t\le s \le T} \|\state{x}_s^{t,x,u}\|^{p} \right] \le K_p (1+\|x\|^{p} ), \ \forall u \in \calU^{\sf s}[t,T].
		\end{equation}
	\hfill$\lhd$
	\end{lemma}

	By applying H\"{o}lder's inequality, we obtain the estimate for the first order moment, that is, \eqref{ineq:moment_upper} also holds for $p=1$.
	
	The estimate~\eqref{ineq:moment_upper} implies $\bbE[\|\state{x}_T^{t,x,u}\|^p] < +\infty$ for any $p\ge 2$. Note that
	\[ \bbE\left[ \int_t^T \psi_0 (u_s) ds \right] \le m(T-t) . \]
	Hence, the growth condition in $(A_2)$ ensures that the cost functional 
	$J^{\sf s} (t,x,u)$ has a finite value for any $(t,x,u)\in [0,T] \times \bbR^n \times \calU^{\sf s}[t,T]$.}



\section{Continuity of $ H^{\sf s} $}\label{app:continuity_H}
	\begin{lemma}\label{lem:continuity_Hs}
		If $f$ and $\sigma$ satisfy \eqref{eq:Lip_fsigma}, 
		then $H^{\sf s}$ defined by~\eqref{eq:H_s} is continuous on $\bbR^n \times \bbR^n \times \calS^n$.
	\end{lemma}
	\begin{pf}
		Fix $\varepsilon > 0$ and $(x,p,M)\in \bbR^n \times \bbR^n \times \calS^n$.
		By definition of $H^{\sf s}$, there exists $\bar{u}\in \bbU$ such that
		\begin{equation}
		H^{\sf s}(x,p,M) - \varepsilon < -f(x,\bar{u})\cdot p - \frac{1}{2} \mathrm{tr}(\sigma \sigma^\top (x,\bar{u})M) - \psi_0(\bar{u}) .
		\end{equation}
		Therefore, for any $(y,q,N)\in \bbR^n \times \bbR^n \times \calS^n$,
		\begin{align*}
		& H^{\sf s}(x,p,M) - H^{\sf s}(y,q,N) \\
		&\le  -f(x,\bar{u}) \cdot p - \frac{1}{2} \mathrm{tr}(\sigma \sigma^\top (x,\bar{u})M) - \psi_0(\bar{u}) + \varepsilon \\
		&\quad	+f(y,\bar{u}) \cdot q + \frac{1}{2} \mathrm{tr}(\sigma \sigma^\top (y,\bar{u})N) + \psi_0(\bar{u}) \\
		&= f(y,\bar{u})\cdot q - f(x,\bar{u})\cdot p \\
		&\quad + \frac{1}{2} \mathrm{tr}(\sigma \sigma^\top (y,\bar{u})N - \sigma \sigma^\top (x,\bar{u})M) + \varepsilon.
		\end{align*}
		Note that $f$ and $\sigma$ are continuous by \eqref{eq:Lip_fsigma}, and thus there exists $\delta > 0$ such that, for any $(y,q,N)\in B((x,p,M),\delta)$,
		\[ f(y,\bar{u})\cdot q - f(x,\bar{u})\cdot p + \frac{1}{2} \mathrm{tr}(\sigma \sigma^\top (y,\bar{u})N - \sigma \sigma^\top (x,\bar{u})M) < \varepsilon. \]
		Hence, for any $(y,q,N)\in B((x,p,M),\delta)$,
		\[  H^{\sf s}(x,p,M) - H^{\sf s}(y,q,N) < 2\varepsilon. \]
		Similarly, $H^{\sf s}(y,q,N) - H^{\sf s}(x,p,M) < 2\varepsilon$ also holds. This shows the continuity of $H^{\sf s}$.
		\hfill$\Box$
\end{pf}

\ito{
	\section{Viscosity solution}\label{app:viscosity}
	\modi{
	Here, we briefly introduce a viscosity solution~\cite{Nisio2014}.
	Let $H: \bbR^n  \times \bbR^n \times \calS^n\to\bbR$ be a continuous function that satisfies the following condition:
	\begin{equation}\label{ineq:ellipticity}
	  H(x,p,M) \le H(x,p,N), \ \textrm{if} \ M-N \in \calS_+^n .
	\end{equation}
	Consider a second-order partial differential equation
	\begin{align}
	\begin{cases}\label{eq:second_order_pde}
	 -v_t(\tx) + H (x, D_x v(\tx), D_x^2 v(\tx)) = 0, \\
	  \hspace{3cm}(\tx) \in [0,T)\times \bbR^n, \\
	   v(T,x) = g(x), \ x\in \bbR^n .
	 \end{cases}
	\end{align}
	A function $v\in C([0,T] \times \bbR^n)$ is said to be a {\em viscosity subsolution} of \eqref{eq:second_order_pde} if
	\[
	v(T,x) \le g(x), \ \forall x\in \bbR^n
	\]
	and, for any $\phi \in C^{1,2} ([0,T)\times \bbR^n) \cap C([0,T]\times \bbR^n)$,
	\begin{equation}\label{ineq:subsolution}
	 -\phi_t(t_0,x_0) + H(x_0,D_x \phi(t_0,x_0), D_x^2 \phi(t_0,x_0))\leq0
	\end{equation}
	at any global maximum point $(t_0,x_0)\in [0,T)\times \bbR^n$ of $v-\phi$. 
	Similarly,
	a function $v\in C([0,T]\times \bbR^n)$ is said to be a {\em viscosity supersolution} of \eqref{eq:second_order_pde}
	if
	\[
		v(T,x) \ge g(x), \ \forall x\in \bbR^n
	\]
	and, for any $\phi\in C^{1,2}([0,T) \times \bbR^n) \cap C([0,T]\times \bbR^n)$,
	\begin{equation}\label{ineq:supersolution}
	-\phi_t(t_0, x_0) + H(x_0,D_x \phi(t_0, x_0), D_x^2 \phi(t_0, x_0)) \geq0
	\end{equation}
	at any global minimum point $(t_0, x_0)\in [0,T)\times \bbR^n$ of $v-\phi$. 
	Finally, $v$ is said to be a {\em viscosity solution} of \eqref{eq:second_order_pde},
	if it is simultaneously a viscosity subsolution and supersolution.
}
	
	Next, we define the {\em second-order right parabolic superdifferential} and {\em subdifferential}, which \rev{are} used in Lemma~\ref{thm:opt_ctrl_stoc} and Theorem~\ref{thm:feedback}. 
	For $v\in C([0,T] \times \bbR^{\mnr{n}})$ with $T>0$, the {\em second-order right parabolic superdifferential} of $v$ at $(\tx)\in [0,T) \times \bbR^n$ is defined by
	\begin{align*}
	&\mnr{D_{t+,x}^{1,2,+}} v(\tx) \triangleq \Bigl\{ (q,p,M) \in \bbR\times \bbR^n \times \calS^n :\\ &\quad \limsup_{\substack{s\searrow t, s\in[0,T)\\y\rightarrow x} }  \frac{1}{|s-t|+\|y-x\|^2} \bigl( v(s,y) - v(t,x)\\
	&\quad -q(s-t) - p\cdot (y-x) - \frac{1}{2} (y-x)^\top M (y-x) \bigr) \le 0    \Bigr\} .
	\end{align*}
	Similarly, the {\em second-order right parabolic subdifferential} of $v$ at $(\tx)\in [0,T)\times \bbR^n$ is defined by
	\begin{align*}
	&D_{t+,x}^{1,2,-} v(\tx) \triangleq \Bigl\{ (q,p,M) \in \bbR\times \bbR^n \times \calS^n :\\ &\quad \liminf_{\substack{s\searrow t, s\in[0,T)\\y\rightarrow x} }  \frac{1}{|s-t|+\|y-x\|^2} \bigl( v(s,y) - v(t,x)\\
	&\quad -q(s-t) - p\cdot (y-x) - \frac{1}{2} (y-x)^\top M (y-x) \bigr) \ge 0    \Bigr\} .
	\end{align*}
	If $v$ admits $v_t, D_x v$, and $D_x^2 v$ at $(t_0,x_0) \in (0,T)\times \bbR^n$,
	it holds that
	\begin{align}
	&\left(v_t (t_0,x_0), D_x v(t_0,x_0), D_x^2 v (t_0,x_0)\right) \nonumber\\
	&\qquad \in D_{t+,x}^{1,2,+} v(t_0,x_0) \cap D_{t+,x}^{1,2,-} v(t_0,x_0) . \label{eq:subdiff_property}
	\end{align}
}


\vspace{-.8cm}

\bibliographystyle{plain}        
\bibliography{auto_stoc_L0}           



\end{document}